\documentclass{amsart}
\usepackage{amssymb}\usepackage{amsmath}\usepackage{a4}
\newtheorem{theorem}{Theorem}[section]
\newtheorem{lemma}[theorem]{Lemma}
\newtheorem{remark}[theorem]{Remark}

\newtheorem{corollary}[theorem]{Corollary}

\begin{document}
\makeatletter
  \renewcommand{\theequation}{%
   \thesection.\alph{equation}}
  \@addtoreset{equation}{section}
 \makeatother
\title[Finite Fourier series]
{Spectral geometry, homogeneous spaces, and differential forms with finite Fourier series}
\author{C. Dunn, P. Gilkey, and J.H. Park}

\address{CD:Mathematics Department, California State University at San Bernardino,
San Bernardino, CA 92407, USA}
\email{cmdunn@csusb.edu}
\address{PG:Mathematics Department, University of Oregon, Eugene, OR 97403, USA}
\email{gilkey@uoregon.edu}
\address{J-H.P:Department of Mathematics, SungKyunKwan University,
Suwon, 440-746, SOUTH KOREA}
\email{parkj@skku.edu}
\begin{abstract} Let $G$ be a compact Lie group acting transitively on Riemannian
manifolds $M_i$ and let $\pi:M_1\rightarrow M_2$ be a $G$-equivariant Riemannian submersion. We show that
a smooth differential form $\phi$ on $M_2$ has finite Fourier series on $M_2$ if and only if the pull-back
$\pi^*\phi$ has finite Fourier series on $M_1$.
\end{abstract}
\maketitle

\section{Introduction} The spectral geometry of Riemannian submersions has been discussed by
many authors; we refer, for example, to \cite{GLP99} for a more extensive discussion. In
particular, it plays an important role in the study of non-bijective canonical transformations;
see, for example, the discussion in
\cite{LK88}.

Let $M$ be a
compact smooth closed Riemannian manifold of dimension $m$, and let $\Delta_M^p$ be the
Laplace-Beltrami operator acting on the space $C^\infty(\Lambda^pM)$ of smooth $p$-forms.
Let $\operatorname{Spec}(\Delta_M^p)$ be the spectrum of $\Delta_M^p$; this is a discrete
countable set of non-negative real numbers. The associated
eigenspaces
$E(\lambda,\Delta_M^p)$ are finite dimensional and there is a complete orthonormal
decomposition
\begin{equation}\label{eqn-1.a}
L^2(\Lambda^pM)=\oplus_{\lambda\in\operatorname{Spec}(\Delta_M^p)}
E(\lambda,\Delta_M^p)
\end{equation}
which we may use to decompose a smooth $p$-form $\phi$ on $M$ in the form $\phi=\sum_\lambda\phi_\lambda$
where $\phi_\lambda\in E(\lambda,\Delta_M^p)$. We say $\phi$ has {\it finite Fourier series} if
this is a finite sum. If $p=0$ and if $M=S^1$, then
this yields, modulo a slight change of notation, the classical Fourier series decomposition
$f(\theta)=\sum_na_ne^{in\theta}$ and a function has a finite Fourier series in this setting if and
only if it is a trigonometric polynomial. There is an extensive literature on the subject, a few
representative items being
\cite{B97,FS06}.

We say that $M$ is a {\it homogeneous space} if there is a compact Lie group $G$ which
acts transitively on $M$ by isometries; if $H$ is the isotropy subgroup associated
to some point $P\in M$, then we may identify $M=G/H$. We may choose a left-invariant
metric $\tilde g$ on $G$ so $g$ is the induced metric or, equivalently, that $\pi:(G,\tilde
g)\rightarrow(M,g)$ is a Riemannian submersion. The following is
the main result of this paper:
\begin{theorem}\label{thm-1.1}
Let $\pi:G\rightarrow G/H$ where  $H$ is a Lie subgroup of a compact Lie group $G$. Let $\tilde g$
be a left-invariant Riemannian metric on
$G$ and let $g$ be the induced Riemannian metric on $G/H$. Then a $p$-form $\phi$ on $G/H$ has
finite Fourier series on $G/H$ if and only if $\pi^*\phi$ has finite Fourier series on $G$.
\end{theorem}

There is an associated Corollary which is useful in applications. 
\begin{corollary}\label{cor-1.2} Let $G$ be a compact Lie group acting transitively on Riemannian
manifolds $M_1$ and $M_2$. Let $\pi:M_1\rightarrow M_2$ be a $G$-equivariant Riemannian submersion. If
$\phi$ is a smooth $p$-form on $M_2$, then $\phi$ has finite Fourier series on
$M_2$ if and only if $\pi^*\phi$ has finite Fourier series on $M_1$.
\end{corollary}

\begin{remark}\label{rmk-1.3}
\rm The Hopf fibration $\pi:S^{2n+1}\rightarrow\mathbb{CP}^n$ is a
$U(n+1)$ equivariant Riemannian submersion which is an important non-canonical transformation used
to study the Coulumb problem, see, for example, the discussion in \cite{BX}. Corollary
\ref{cor-1.2} shows $\phi$ has finite Fourier series on $\mathbb{CP}^n$ if and only if $\pi^*\phi$
has finite Fourier series on $S^{2n+1}$.
\end{remark}

\section{The proof of Theorem \ref{thm-1.1}}

The central ingredient is our discussion is the classical Peter--Weyl theorem \cite{H03}. Let
$\operatorname{Irr}(G)$ be the collection of equivalence classes of irreducible finite dimensional
representations of
$G$; if $\rho\in\operatorname{Irr}(G)$, let $V_\rho$ be the associated representation space. The Hilbert
space structure on $L^2(G)$ depends on the particular Riemannian metric which is chosen; this space is
invariantly defined as a Banach space, however. This is a minor distinction which will be useful,
however, in Section \ref{sect-4}. Left multiplication defines an action of
$G$ on $L^2(G)$. This action decomposes as a direct sum
\begin{equation}\label{eqn-2.a}
L^2(\Lambda^pG)=\oplus_{\rho\in\operatorname{Irr}(G)}W_\rho
\end{equation}
where each $W_\rho$ is a finite dimensional irreducible subspace of $L^2(G)$ which is isomorphic to a
finite number of copies of
$V_\rho$. If
$\Phi$ is a smooth $p$-form on $G$, we may use Equation (\ref{eqn-2.a}) to decompose
$\Phi=\sum_\rho\Phi_\rho$ for $\Phi_\rho\in W_\rho$. We say that $\Phi$ has {\it finite
representation expansion} on $G$ if this sum is finite; we emphasize that this notion is independent of
the particular Riemannian metric chosen.

Since $\pi$ is a submersion, $\pi^*$ is an injective $G$-equivariant map from
$L^2(\Lambda^p(G/H))$ to $L^2(G)$ with closed image. The decomposition
$$L^2(\Lambda^pG)=\pi^*(L^2(\Lambda^p(G/H)))\oplus\{\pi^*(L^2(\Lambda^p(G/H)))\}^\perp$$
is $G$-equivariant. We therefore have an orthogonal direct sum decomposition of
$L^2(\Lambda^p(G/H))$ as a representation space for $G$ in the form:
\begin{eqnarray}\label{eqn-2.b}
&&L^2(\Lambda^p(G/H))=\oplus_{\rho\in\operatorname{Irr}(G)}X_\rho\quad\hbox{where}\\
&&\pi^*X_\rho=W_\rho\cap\pi^*(L^2(\Lambda^p(G/H)))\,.\label{eqn-2.c}
\end{eqnarray}
We say that a $p$-form $\phi$ on $G/H$ has {\it finite $G$-representation series} if the expansion
$\phi=\sum_\rho\phi_\rho$ given by Equation (\ref{eqn-2.b}) is finite. Theorem \ref{thm-1.1} will
follow from the following:

\begin{lemma}\label{lem-2.1}
Adopt the notation established above. Let $\phi$ be a smooth $p$-form on $G/H$. Fix a left-invariant
$\tilde g$ metric on $G$ and let $g$ be the induced metric on $G/H$. The following assertions are
equivalent:
\begin{enumerate}
\item $\phi$ has finite Fourier series on $G/H$.
\item $\phi$ has finite $G$-representation series on $G/H$.
\item $\pi^*\phi$ has finite Fourier series on $G$.
\item $\pi^*\phi$ has finite $G$-representation series on $G$.
\end{enumerate}\end{lemma}

\begin{proof} The equivalence of Assertions (ii) and (iv) is immediate from
Equation (\ref{eqn-2.c}). We argue as follows to prove that Assertion (i) implies Assertion (ii).
Suppose that $\phi$ has finite Fourier series on $G/H$. Since $G$ acts by isometries, $G$ commutes
with the Laplacian. Thus $E(\lambda,\Delta_{G/H}^p)$ is a finite dimensional representation space
for
$G$. Only a finite number of representations occur in the representation decomposition of
$E(\lambda,\Delta_{G/H}^p)$ and thus any eigen
$p$-form on
$G/H$ has finite
$G$-representation series on
$G/H$; more generally, of course, any finite sum of eigen $p$-forms on $G/H$ has finite
$G$-representation series on $G/H$. This shows that Assertion (i) implies Assertion (ii); a similar
argument shows Assertion (iii) implies Assertion (iv).

Each representation appears with finite multiplicity in $L^2(\Lambda^p(G/H))$. Thus each
representation appears in the decomposition of $E(\lambda,\Delta_{G/H}^p)$ for only a finite
number of $\lambda$. Thus any element of $X_\rho$ has finite Fourier series and more generally any
$p$-form on $G/H$ with finite $G$-representation series has finite Fourier series. Thus
Assertion (ii) implies Assertion (i); similarly, Assertion (iv) implies Assertion (iii).
\end{proof}
\section{The proof of Corollary \ref{cor-1.2}}
Let $\pi:M_1\rightarrow M_2$ be a $G$-equivariant Riemannian submersion; this means that we may
express               
$M_i=G/H_i$ where
$H_1\subset H_2\subset G$. Let $\pi_i:G\rightarrow G/H_i$ be the natural projections. We then have
$\pi\pi_1=\pi_2$ and thus $\pi_2^*=\pi_1^*\pi^*$. Let $\phi$ be a smooth $p$-form on $G/H_2$. We
apply Theorem
\ref{thm-1.1} to derive the following chain of equivalent statements from which Corollary
\ref{cor-1.2} will follow:
\begin{enumerate}
\item $\phi$ has finite Fourier series on $G/H_2$.
\item $\pi_2^*\phi$ has finite Fourier series on $G$.
\item $\pi_1^*(\pi^*\phi)$ has finite Fourier series on $G$.
\item $\pi^*\phi$ has finite Fourier series on $G/H_1$.
\end{enumerate}
\section{Conclusions and open problems}\label{sect-4}
Our methods in fact show a bit more. Let $g_i$ be two left invariant metrics on $G$ and let $\phi$
be a smooth $p$-form on $G$. Then $\phi$ has finite Fourier series with respect to $g_1$ if and
only if $\phi$ has finite Fourier series with respect to $g_2$ since both conditions are
equivalent to $\phi$ having finite representation series and this notion is independent of the
particular metric chosen.

Cayley multiplication defines a Riemannian submersion $\pi:S^7\times S^7\rightarrow S^7$. The
group of isometries commuting with this action does not, however, act transitively on $S^7\times
S^7$ and Theorem \ref{thm-1.1} is not applicable. Our
research continues in this area as this example has important physical applications (see, for
example,
\cite{LK88}).
\section*{Acknowledgments}
Research of C. Dunn partially supported by a CSUSB faculty research
grant. Research of P. Gilkey partially supported by the Max Planck
Institute in the Mathematical Sciences (Leipzig, Germany). Research of
both C. Dunn and P. Gilkey partially supported by the University of Santiago (Spain). Research of J.H.Park
partially supported by the Korea Research Foundation Grant funded by the Korean Government (MOEHRD)
KRF-2007-531-C00008

\end{document}